\documentclass[11pt]{amsart}
\usepackage{amssymb,diagrams}

\newcommand{\OnS}{Ozsv\'ath and Szab\'o}

\newcommand{\zee}{\mbox{$\mathbb Z$}}
\newcommand{\cee}{\mbox{$\mathbb C$}}
\newcommand{\cue}{\mbox{$\mathbb Q$}}
\newcommand{\eff}{\mbox{$\mathbb F$}}
\newcommand{\spinc}{\mbox{{spin}${}^c$ }} 
\newcommand{\s}{\mathfrak s}
\newcommand{\R}{\mathfrak r}

\newcommand{\rk}{{\mbox{rk}}}

\newcommand{\cH}{{\mbox{$\mathcal H$}}}

\newcommand{\Tor}{{\mbox{\rm Tor}}}
\newcommand{\E}{{\mathcal E}}

\newcommand{\ds}{\displaystyle}

\newcommand{\ts}{\textstyle}
\newcommand{\longto}{\longrightarrow}
\newcommand{\tot}{\mbox{tot}}
\newtheorem{thm}{Theorem}
\newtheorem{cor}[thm]{Corollary}
\newtheorem{lemma}[thm]{Lemma}

\newtheorem{prop}[thm]{Proposition}
\newtheorem{defn}[thm]{Definition}
\newtheorem{remark}[thm]{Remark}

\begin{document}

\title[Triple products and invariants for three-manifolds]{Triple products and cohomological invariants for closed three-manifolds}
\author{Thomas E. Mark}

\begin{abstract} Motivated by conjectures in Heegaard Floer homology, we introduce an invariant $HC^\infty_*(Y)$ of the cohomology ring of a closed 3-manifold $Y$ whose behavior mimics that of the Heegaard Floer homology $HF^\infty(Y,\s)$ for $\s$ a torsion \spinc structure. We derive from this a numerical invariant $h(Y)\in\zee$, and obtain upper and lower bounds on $h(Y)$. We describe the behavior of $h(Y)$ under connected sum, and deduce some topological consequences. Examples show that the structure of $HC^\infty_*(Y)$ can be surprisingly complicated, even for 3-manifolds with comparatively simple cohomology rings.
\end{abstract}
\maketitle

Heegaard Floer homology groups are a powerful tool in 
low-dimensional topology introduced and studied by {\OnS} 
(\cite{OSplumb}, \cite{OS1}, \cite{OS2}, etc.), and they have generated much interest among topologists (\cite{hedden}, \cite{stipsicz}, \cite{nemethi}, etc.). The groups are associated to a closed oriented 
3-manifold $Y$ together with a choice of \spinc structure $\s$, and 
comprise a number of variations: $HF^+$, $HF^-$, $\widehat{HF}$, $HF^\infty$. Of these, $HF^\infty$ is considered to be the least interesting as an invariant, due to the apparent fact (formulated as a conjecture by {\OnS} \cite{OSplumb}) that it is determined by the cohomology ring of $Y$. However, while all evidence supports Ozsv\'ath's and Szab\'o's conjecture, the structure of $HF^\infty$ can be rather more complicated than a cursory inspection of the cohomology ring of $Y$ might suggest (c.f. \cite{sigxs1}). Furthermore, in various situations it can be useful for other purposes to understand the behavior of $HF^\infty$---for example, it plays a key role in Ozsv\'ath's and Szab\'o's proof of Donaldson's diagonalizability theorem for definite 4-manifolds and generalizations \cite{OS4}.

With these ideas in mind, we introduce here an invariant $HC^\infty_*(Y)$ of the cohomology ring of $Y$ that we call the ``cup homology;'' it is closely related to $HF^\infty(Y,\s)$ for any torsion \spinc structure $\s$, granted the conjecture mentioned above (more precisely, in this case $HC^\infty_*(Y)$ is the $E^\infty$ term of a spectral sequence converging to $HF^\infty(Y,\s)$, possibly after a grading shift---see section \ref{floersec} below). The cup homology satisfies various properties:
\begin{enumerate}
\item It is the homology of a free complex $C^\infty_*(Y)$ over $\zee$ whose underlying group is $\Lambda^*H^1(Y;\zee)\otimes \zee[U,U^{-1}]$ (where $U$ is a formal variable of degree $-2$), and whose differential is defined in terms of triple products $\langle a\cup b\cup c,[Y]\rangle$ of elements $a,b,c\in H^1(Y)$.
\item Multiplication by $U$ induces an isomorphism between $HC^\infty_k(Y)$ and $HC^\infty_{k-2}(Y)$ for all $k$, so $HC^\infty_*(Y)$ is determined as a $\zee[U]$-module by its values in two adjacent degrees.
\item When $b_1(Y)\geq 1$, the ranks $\rk_{\mathbb Z}(HC^\infty_{2k}(Y))$ and $\rk_{\mathbb Z}(HC^\infty_{2k+1}(Y))$ are equal.
\end{enumerate}
The reader familiar with Heegaard Floer homology will recognize that $C^\infty_*(Y)$ is identical (at least on the level of groups) with the $E_2$ term of a spectral sequence that calculates $HF^\infty(Y,\s)$. (The $d_2$ differential vanishes for trivial reasons, and according to the conjecture of \OnS, $d_3$ is given by the boundary operator used here.)

We define an invariant $h(Y)$ of the cohomology ring of $Y$ by:
\[
h(Y) = \left\{\begin{array}{ll} \rk_{\mathbb Z} HC^\infty_k(Y), &\mbox{any $k\in\zee$, for $b_1(Y) \geq 1$}\\
\frac{1}{2} & b_1(Y) = 0.\end{array}\right.
\]

Our goal here is to investigate the possibilities for $HC^\infty_*(Y)$, and in particular to study the behavior of $h(Y)$. 

Suppose $b_1(Y)\leq 2$. Then there can be no nontrivial triple products of elements in $H^1(Y;\zee)$, so the differential on $C^\infty_*(Y)$ vanishes. It follows that in this case the group $HC^\infty_*(Y)$ is independent of $Y$, and
\begin{eqnarray*}
\mbox{If $b_1(Y) = 0$}:&&HC^\infty_{2k}(Y) = \zee, \quad HC^\infty_{2k+1}(Y) = 0\\
\mbox{If $b_1(Y) = 1$}:&& HC^\infty_k(Y) = \zee\mbox{ for all $k$, so } h(Y) = 1\\
\mbox{If $b_1(Y) = 2$}:&& HC^\infty_k(Y) = \zee^2\mbox{ for all $k$, so } h(Y) = 2.
\end{eqnarray*}

In general, we have the following bounds on $h(Y)$.

\begin{thm} Fix an integer $b\geq 1$. Then for any 3-manifold 
$Y$ with $b_1(Y) = b$ we have
\begin{equation}
L(b) \leq h(Y) \leq 2^{b-1}
\label{mainbound}
\end{equation}
where
\[
L(b) = \left\{\begin{array}{ll} 3^{\frac{b-1}{2}} & \mbox{if $b$ 
is odd}\\ 2\cdot 3^{\frac{b}{2} - 1} & \mbox{if $b$ is 
even}.\end{array}\right.
\]
\label{mainthm}
\end{thm}

It is natural to ask, given a fixed value $b$ of $b_1(Y)$, which values of $h(Y)$ allowed by 
Theorem \ref{mainthm} can be realized. We can certainly realize the upper bound $2^{b-1}$ 
for any $b\geq 1$ by taking $Y$ to be the connected sum of $b$ copies of $S^1\times 
S^2$, or any other $3$-manifold with $b_1(Y) = b$ having trivial cup products on $H^1$. In the case $b = 3$, 
the two possibilities $h = 3$ and $h = 4$ are realized by the 
3-torus $T^3$ and $\#^3S^1\times S^2$, respectively. On the other hand, a simple calculation based on the definitions (Lemma \ref{constrlemma}) shows that when $b_1(Y) \geq 4$, the value $h = 2^{b_1(Y) - 1} -1$ does not occur. Further examples are discussed in section \ref{examplesec} below.
In general, one expects that a more complicated cup product 
structure in the cohomology of $Y$ will result in a smaller value for 
$h(Y)$, where ``complicated'' refers roughly to the number of 
nonvanishing triple products of elements of $H^1(Y;\zee)$. It is a result of D. Sullivan \cite{sullivan} that in an appropriate sense, all possibilities for triple-cup-product behavior are realized by closed 3-manifolds $Y$ (see section \ref{prelimsec} below for a more precise statement). In principle, this result reduces the determination of which values of $h$ are realized to a purely algebraic-combinatorial question.

We have the following general result on the behavior of $h$ under connected sum.

\begin{thm} Let $Y_1$ and $Y_2$ be closed 3-manifolds as above. Then
\[
h(Y_1\# Y_2) = 2h(Y_1)h(Y_2).
\]\label{2ndthm}
\end{thm}

Thus, for example, $h(Y\# S^1\times S^2) = 2h(Y)$, which gives an easy way to realize values of $h$ recursively. To spell this out, let $\cH_b\subset \{L(b), L(b) + 1,\ldots, 2^{b-1}\}$ denote the collection of integers $n$ such that $n = h(Y)$ for some 
3-manifold $Y$ with $b_1(Y) = b$. Then:

\begin{cor} For any integer $b\geq 1$, the set $\cH_{b+1}$ contains
$2\cH_b = \{2h | h\in \cH_b\}$. In fact, the operation $Y\mapsto Y\#S^1\times 
S^2$ demonstrates an inclusion
\begin{eqnarray*}
\cH_b&\to&\cH_{b+1}\\
h&\mapsto& 2h.
\end{eqnarray*}
\label{evencor}
\end{cor}

It follows, for example, that if $Y$ is a 3-manifold with $b_1(Y) = 
b$ realizing the lower bound $h(Y) = L(b)$ for $b$ odd, then $Y\# 
S^1\times S^2$ realizes the lower bound $h(Y\# S^1\times S^2) = L(b+1)$. However, it is clear from Theorem \ref{mainthm} that this naive construction will not help realize small values of $h(Y)$.

\begin{cor} The rank of $HF^\infty_k(T^3\#S^1\times S^2; \s_0)$ (for 
$\s_0$ the torsion {\spinc} structure)
realizes the smallest possible rank $L(4)=6$ for $HF^\infty$ in each degree 
among all 3-manifolds $Y$ having $b_1(Y) = 4$.
\label{minhcor}
\end{cor}

\begin{proof} We will see below (section \ref{floersec}) that the lower bound on $h(Y)$ of Theorem \ref{mainthm} also gives a lower bound on the rank of $HF^\infty_k(Y,\s)$ in any torsion \spinc structure $\s$, whenever $b_1(Y)\leq 4$. Hence the corollary follows from Theorems \ref{mainthm} and \ref{2ndthm}.
\end{proof}

Note that the corollary could easily be proved directly in Heegaard Floer homology using a 
spectral sequence argument as in the proof of Lemma 4.8 of 
\cite{OSplumb}. 

A more subtle question in regards to the 3-manifold ``geography'' 
we're considering is: which values of $h(Y)$ can be realized by 
{\it irreducible} 3-mani\-folds? For example, we will see below that if $\Sigma_g$ is a closed orientable surface of genus $g\geq 1$ then 
\[
h(\Sigma_g\times S^1) = {{2g+1}\choose{g}},
\]
in particular $h(\Sigma_2\times S^1) = 10$. Note that since $b_1(\Sigma_2\times S^1) = 5$ the bounds supplied by Theorem \ref{mainthm} are $9\leq h\leq 16$ (the analogue of the latter fact for $HF^\infty$ is part of the content of Lemma 4.8 of \cite{OSplumb}).

Let us call a 3-manifold $Y$ {\it rationally irreducible} if in any connected sum decomposition of $Y$, at least one of the factors is a rational homology sphere. Observe that since connect sum with a rational homology sphere does not change the triple-cup-product structure of a 3-manifold, if an integer $h$ is realized as $h(Y)$ for some rationally irreducible 3-manifold $Y$ then it is also realized by an irreducible 3-manifold.

We have the following immediate corollary of Theorem \ref{2ndthm}:

\begin{cor}\label{oddcor} If $h(Y)$ is odd, then $Y$ is rationally irreducible.
\end{cor}

We also have the following result related to irreducibility, 
which follows from the  
behavior of $h$ under connected sum:

\begin{thm}\label{irredthm}
Let $Y_1$ and $Y_2$ be 3-manifolds having first Betti
number at least 1, and suppose that $b_1(Y_1)$ and $b_1(Y_2)$ are not both
odd. Let $Y = Y_1\# Y_2$, and set $b = b_1(Y) = b_1(Y_1)+b_1(Y_2)$.
Then
\[
h(Y) \geq \frac{4}{3} L(b).
\]
\end{thm}

As a consequence we see that if one is interested in realizing small
values of $h(Y)$ for a fixed, odd, value of $b_1$, then the only
candidates are (rationally) irreducible.

On the other hand, it follows from the theorem that if $Y$ is a 3-manifold with $b=b_1(Y)$ odd and $h(Y) < \frac{4}{3}L(b) = 4\cdot 3^{(b-3)/2}$ then $Y$ is rationally irreducible.
The inequality in Theorem \ref{irredthm} cannot be strengthened, as shown by the
example
\[
T^3\,\#\, S^1\times S^2 \,\#\, S^1\times S^2,
\]
which has $h = 12 = \frac{4}{3}L(5)$.

The complex $C^\infty_*(Y)$ can be just as easily defined using cohomology with coefficients in any commutative ring, e.g., using $H^1(Y;\zee_p)$, giving rise to a mod-$p$ cup homology $HC^\infty_*(Y)_p$. We obtain a sequence of invariants $h_p(Y)$ given by the rank in each dimension of $HC^\infty_*(Y)_p$; the arguments used to define $h(Y)$ and obtain Theorem \ref{mainthm} and the other results listed here are insensitive to this change (where of course $b_1(Y)$ is calculated in the appropriate coefficients). These invariants can easily be distinct from each other and from $h(Y)$ (see section \ref{examplesec} below), and we can consider the realization problem for each of them.

It follows from the results above that if we define for any $p\geq 1$
\[
k_p(Y) = \log_2(2h_p(Y))
\]
(where by convention $h_1(Y) = h(Y)$), then $\{k_1(Y),k_2(Y), k_3(Y),\ldots\}$ is a sequence of real-valued invariants of 3-manifolds that vanish for rational homology spheres, are additive under connected sum, and satisfy 
\[
m(b_1(Y))\leq k_p(Y) \leq b_1(Y)
\]
where the lower bound $m$ is a linear function of the first Betti 
number easily derived from Theorem \ref{mainthm}.

In the next section we define the chain complex $C^\infty_*(Y)$ and make some simple observations; section \ref{examplesec} is devoted to a few sample calculations. We prove Theorem \ref{mainthm} in section \ref{thm1proof}, and in section \ref{connsumproof} we prove Theorems \ref{2ndthm} and \ref{irredthm}. In the final section we spell out the conjectural relationship between $HC^\infty_*(Y)$ and $HF^\infty(Y,\s)$. One remark is in order here, namely that there is another version of Floer homology for 3-manifolds based on the Seiberg-Witten equations, due to Kronheimer and Mrowka \cite{KM}. This ``monopole Floer homology'' appears to be isomorphic to Heegaard Floer homology, and furthermore the relationship between the cohomology ring and the Seiberg-Witten analog of $HF^\infty$ (that is, $\overline{HM}_*$)  has been established in that theory. Therefore, according to chapter IX of \cite{KM}, our $HC^\infty_*$ is isomorphic to the $E_\infty$ term of a spectral sequence converging to the monopole Floer homology $\overline{HM}_*$, modulo a possible shift in grading. In particular, most of our results here could be rephrased using this monopole homology rather than the---more elementary but somewhat artificial---cup homology.

\section{Definitions}\label{prelimsec}
Let $Y$ be a closed oriented 3-dimensional manifold and write $H = H^1(Y;\zee)$. The cup product structure on $Y$ induces a 3-form on $H$, written $\mu_Y\in \Lambda^3H^*$, given by $\mu_Y(a,b,c) = \langle a\cup b\cup c, [Y]\rangle$ for $a,b,c\in H$.

\begin{remark} Sullivan \cite{sullivan} has shown that for any pair $(H,\mu)$ where $H$ is a finitely generated free abelian group and $\mu\in \Lambda^3H^*$, there exists a 3-manifold $Y$ with $H^1(Y) =H$ and cup-product form $\mu_Y = \mu$. 
\end{remark}

 There is a natural interior product $H^*\otimes \Lambda^kH \to \Lambda^{k-1}H$, written $\omega\otimes \alpha\mapsto \omega\angle\alpha$, induced by the duality between $H^*$ and $H$, with the property that $\omega\angle(\omega\angle \alpha) = 0$ for $\omega\in H^*$. Thus there is an extension of $\angle$ to the exterior algebra:
\[
\Lambda^pH^*\otimes \Lambda^kH\to \Lambda^{k-p}H,
\]
satisfying $(\omega\wedge\eta)\angle \alpha = \omega\angle(\eta\angle\alpha)$. 

\begin{defn} Let $U$ be a formal variable of degree $-2$. The {\em cup complex} of $Y$ is defined to be the chain complex $C^\infty_* = C^\infty_*(Y)$ with chain groups
\[
C_*^\infty = \Lambda^*H \otimes \zee[U, U^{-1}]
\]
graded in the obvious way, and with differential
\begin{eqnarray}
&\partial : C_k^\infty\to C_{k-1}^\infty&\nonumber\\
&\partial (\alpha\otimes U^n) = \mu_Y\angle \alpha\otimes U^{n-1},&\label{diffdef}
\end{eqnarray}
where $\mu_Y$ is the 3-form given by cup product defined above.
\end{defn}

Thus for any $k$,
\begin{equation}
C^\infty_k \cong \bigoplus_{\ell\equiv k \mod 2} \Lambda^\ell H \otimes U^{(\ell -k)/2},
\label{chaingp}
\end{equation}
and
\[
\partial: \Lambda^\ell H\otimes U^n\to \Lambda^{\ell-3}\otimes U^{n-1}.
\]
We have an explicit expression for $\partial$, namely if $a_1,\ldots,a_k\in H^1(Y)$, 
\[
\partial(a_1\cdots a_k) = U^{-1} \cdot \sum_{i_1<i_2<i_3} (-1)^{i_1+i_2+i_3} \langle a_{i_1}\cup a_{i_2}\cup a_{i_3},[Y]\rangle a_1\cdots\widehat{a_{i_1}}\cdots\widehat{a_{i_2}}\cdots\widehat{a_{i_3}}\cdots{a_k}
\]
where we use juxtaposition to indicate wedge product.

Since $\mu_Y\angle(\mu_Y\angle \alpha) = (\mu_Y\wedge\mu_Y)\angle \alpha = 0$, we see $(C_*^\infty(Y),\partial)$ is indeed a chain complex. We write $HC_*^\infty(Y)$ for $H_*(C_*^\infty(Y),\partial)$; clearly $HC_*^\infty(Y)$ is an invariant of the homotopy type of $Y$.

Obviously $U: C^\infty_*\to C^\infty_{*-2}$ is a chain isomorphism, so the homology $HC^\infty_*(Y)$ is determined as a group by its values in two adjacent degrees. Let 
\begin{eqnarray*}
&h_{ev}(Y) = \rk_{\mathbb Z}(HC^\infty_{2k}(Y))&\\
&h_{odd}(Y) = \rk_{\mathbb Z}(HC^\infty_{2k+1}(Y))&
\end{eqnarray*}
for any $k$. Observe that $HC^\infty_{ev} \cong HC_{2k}\otimes \zee[U,U^{-1}]$ and $HC^\infty_{odd} \cong HC_{2k+1}(Y) \otimes \zee[U,U^{-1}]$ for any $k$.

\begin{lemma} If $b_1(Y)\geq 1$, then $h_{ev}(Y) = h_{odd}(Y)$.\label{eulerlemma}
\end{lemma}

\begin{proof} As a chain complex over the graded ring $\zee[U,U^{-1}]$, the complex $C^\infty_*$ has even- and odd- graded parts given by $\Lambda^{ev}H\otimes \zee[U,U^{-1}]$ and $\Lambda^{odd}H\otimes \zee[U,U^{-1}]$. When $b_1(Y)\geq 1$ the groups $\Lambda^{ev}H$ and $\Lambda^{odd}H$ have the same rank, so that $C^\infty_*$ has vanishing Euler characteristic over $\zee[U,U^{-1}]$. Therefore the same is true of $HC^\infty_*$, so that $0 = \rk_{{\mathbb Z}[U,U^{-1}]}(HC^\infty_{ev}) - \rk_{{\mathbb Z}[U,U^{-1}]}(HC^\infty_{odd}) = h_{ev} - h_{odd}$.
\end{proof}

\begin{defn} For $Y$ a closed oriented 3-manifold with $b_1(Y)\geq 1$, define
\[
h(Y) = h_{ev}(Y) = h_{odd}(Y).
\]
If $b_1(Y) = 0$, set $h(Y) = \frac{1}{2}$.
\end{defn}

Thus $h$ is an invariant of the cohomology ring of $Y$ that takes values in $\zee$ when $Y$ is not a rational homology sphere.

\begin{lemma}\label{ublemma} For any 3-manifold $Y$, we have $h(Y)\leq 2^{b_1(Y)-1}$.
\end{lemma}

\begin{proof} The lemma is true by definition if $b_1(Y) = 0$, so suppose $b_1(Y)\geq 1$. The rank of $HC^\infty_k(Y)$ is no larger than the rank of $C^\infty_k$, so from \eqref{chaingp} we have
\[
\rk( HC^*_k(Y)) \leq\sum_{\ell\equiv k \mod 2} \rk(\Lambda^\ell H) = \left\{\begin{array}{ll} \ds\sum_{\ell\mbox{\scriptsize { even}}} {{b_1(Y)}\choose{\ell}} & \mbox{$k$ even}\vspace{1em}\\
\ds\sum_{\ell\mbox{\scriptsize { odd}}} {{b_1(Y)}\choose{\ell}} & \mbox{$k$ odd}\end{array}\right.
\]
Since both sums on the right are equal to $2^{b_1(Y)-1}$, the lemma follows.
\end{proof}

Following similar terminology in Heegaard Floer homology, we say that a 3-manifold $Y$ has {\it standard} cup homology if $h(Y) = 2^{b_1(Y) - 1}$. Equivalently, the cup homology is standard if and only if the triple product form $\mu_Y$ is zero.

The following is another example of a constraint on $h$ arising from purely algebraic considerations.

\begin{lemma}\label{constrlemma} If $Y$ is a 3-manifold with $b_1(Y) = b$, $b\geq 4$, and $h(Y)\neq 2^{b-1}$ then $h(Y)\leq 2^{b-1}-2$.
\end{lemma}

\begin{proof} Suppose the triple product form on $Y$ is given by $\mu_Y = \sum a_{ijk} e_ie_je_k \in \Lambda^3H^*$, for a basis $\{e_n\}$ of $H^*$ and $a_{ijk}$ integers, with the sum over $0<i<j<k\leq b$. Then the differential acts on the top exterior power $\Lambda^bH \otimes U^n$ by
\[
\partial : e_1\cdots e_b\otimes U^n \mapsto \sum_{i<j<k} \pm a_{ijk} e_1\cdots \hat{e}_i\cdots \hat{e}_j\cdots \hat{e}_k\cdots e_b\otimes U^{n-1},
\]
which is injective unless $\mu_Y = 0$. Likewise, the component of the differential mapping into $\Lambda^0H$ is easily seen to be surjective (over the rationals) if and only if $\mu_Y\neq 0$. 

Therefore as soon as $\mu_Y\neq 0$ a 2-dimensional space of cycles disappears from $\Lambda^bH\oplus \Lambda^3H$ and a 2-dimensonal space of boundaries appear in $\Lambda^{b-3}H\oplus \Lambda^0H$. This implies that the rank of $HC^\infty_{2k}(Y)$ and $HC^\infty_{2k+1}(Y)$ are each forced to decrease by at least 2 once $b>3$ and $\mu_Y \neq 0$.
\end{proof}

\section{Examples}\label{examplesec}

Here we calculate $h(Y)$ for some sample 3-manifolds $Y$. As remarked previously, since there can be no nontrivial triple products on $H^1(Y)$ when $b_1(Y)\leq 2$, the differential on $C^\infty_*(Y)$ must vanish in this case. This easily gives the results listed before Theorem \ref{mainthm}.

More generally, if all triple cup products of elements of $H^1(Y;\zee)$ vanish---i.e., if $\mu_Y = 0$---then the cup homology is standard, i.e., $HC^\infty_*(Y) \cong \Lambda^*H^1(Y)\otimes\zee[U,U^{-1}]$. This is the case, for example, for connected sums of copies of $S^1\times S^2$. For more interesting examples, it was observed by Sullivan \cite{sullivan} that if $Y$ is the link of an isolated algebraic surface singularity then $\mu_Y = 0$. In fact, Sullivan's argument proves:

\begin{prop} If $Y$ is a closed 3-manifold bounding an oriented 4-man\-i\-fold $X$ such that the cup product pairing on $H^2(X,Y)$ is nondegenerate, then $\mu_Y = 0$. Hence $HC^\infty_*(Y) $ is standard for such $Y$, so that $h(Y) = 2^{b_1(Y)-1}$.
\end{prop}

In particular, since the link of a singularity bounds a 4-manifold with negative-definite intersection form on $H_2(X) = H^2(X,Y)$, such $3$-manifolds have standard $HC^\infty_*(Y)$.

\begin{proof} We work over the rationals, since torsion cannot contribute to triple products over $\zee$. Poincar\'e duality implies that the cup pairing
\[
Q_X : H^2(X)\otimes H^2(X,Y)\to \cue
\]
is nondegenerate, i.e., there is an isomorphism $Q_X: H^2(X)\to (H^2(X,Y))^*$. Let $\tilde{Q}_X: H^2(X,Y)\otimes H^2(X,Y)\to \cue$ be the cup product form on $H^2(X,Y)$; then we have a commutative diagram
\begin{diagram}
&&H^2(X)&&\\
&\ruTo^{i} & & \rdTo^{Q_X}_\sim & \\
H^2(X,Y) & &\rTo_{\tilde{Q}_X} && (H^2(X,Y))^*
\end{diagram}
Hence $\tilde{Q}_X$ is nondegenerate if and only if $i$ is an isomorphism. In this case, in the sequence
\[
\cdots\to H^1(Y)\stackrel{j}{\to} H^2(X,Y)\stackrel{i}{\to} H^2(X)\stackrel{k}{\to} H^2(Y)\to\cdots
\]
the homomorphisms $j$ and $k$ vanish. As noted in \cite{sullivan}, this means that cup products of elements $a,b\in H^1(Y)$ can be computed by lifting to $H^1(X)$ since $j$ is trivial, multiplying, and restricting to $Y$, which gives 0 since $k$ is trivial.

\end{proof}

If $b_1(Y) = 3$ the only possible nontrivial differential in $C^\infty_*(Y)$ is between $\Lambda^3H\otimes U^n$ and $\Lambda^0H\otimes U^{n-1}$. If $a\cup b\cup c = 0$ for a basis $a,b,c$ of $H = H^1(Y)$, then we have $HC^\infty_k(Y) \cong \zee^4$ for each $k$, so $h(Y) = 4$. Otherwise, i.e., if $a\cup b\cup c \neq 0$, the differential is injective and $HC^\infty_k(Y)$ has rank 3 for all $k$, hence $h(Y) = 3$. In fact, if $\langle a\cup b\cup c, [Y]\rangle = n\neq 0$ then
\begin{eqnarray*}
HC^\infty_{2k}(Y) &\cong& \zee^3 \oplus (\zee/n\zee)\\
HC^\infty_{2k+1}(Y) &\cong& \zee^3.
\end{eqnarray*}
It also follows that in the notation of the introduction, $h_p(Y) = 4$ for primes $p$ dividing $n$, while $h_q(Y) = 3$ for other primes $q$.

A more substantial example is given by $Y = \Sigma_g\times S^1$, where $\Sigma_g$ is a closed oriented surface of genus $g\geq 1$. In this case it is a simple matter to calculate that $\mu_Y = s\wedge \omega$, where $s$ is the class $[pt\times S^1]$ and $\omega\in \Lambda^2H_1(\Sigma_g)$ is the symplectic 2-form given by the cup-product pairing on the first cohomology of $\Sigma$. We have a decomposition
\[
\Lambda^k H^1(\Sigma_g\times S^1) \cong \Lambda^kH^1(\Sigma_g) \oplus \left(\sigma\wedge \Lambda^{k-1}H^1(\Sigma_g)\right),
\]
where $\sigma$ is Poincar\'e dual to $[\Sigma_g\times pt]$, and with respect to this decomposition, the boundary in $C^\infty_*(\Sigma_g\times S^1)$ is trivial on the first factor and given by
\[
\partial(\sigma\wedge\alpha\otimes U^n) = \omega\angle \alpha\otimes U^{n-1} \in \Lambda^{k-2}H^1(\Sigma_g)\otimes U^{n-1}
\]
on the second factor, for $\alpha\in \Lambda^kH^1(\Sigma_g)$. Let $E_k(g)$ denote the Abelian group determined by the long exact sequence
\[
\cdots \Lambda^{k+1}H^1(\Sigma_g)\stackrel{\omega\angle\cdot}{\to} \Lambda^{k-1}H^1(\Sigma_g)\to E_k(g) \to \Lambda^{k}H^1(\Sigma_g)\stackrel{\omega\angle\cdot}{\to} \Lambda^{k-2}H^1(\Sigma_g)\cdots
\]
That is, with appropriate grading conventions $E_*(g)$ is the homology of the mapping cone of $\omega\angle \cdot$ acting on $\Lambda^*H^1(\Sigma_g)$, thought of as a complex with trivial differential. (Note that since the latter is free Abelian, $E_*(g)$ is uniquely determined by the sequence given.) It follows from the discussion above that
\begin{equation}
HC^\infty_*(\Sigma_g\times S^1) \cong  E_*(g)\otimes \zee[U,U^{-1}].
\label{sigxs1homol}
\end{equation}

To give an explicit expression for $E_k(g)$, observe that there is a natural duality isomorphism $\star:\Lambda^kH^1(\Sigma_g)\to \Lambda^{2g-k}H_1(\Sigma_g)$ induced by interior product with the orientation form $\frac{1}{g!}\omega^g\in \Lambda^{2g}H_1(\Sigma_g)$. Identifying $H_1(\Sigma_g)$ with $H^1(\Sigma_g)$ via Poincar\'e duality, it is an exercise to see that $\star(\omega \angle \alpha) = \omega\wedge\star\alpha$ for any $\alpha\in \Lambda^*H^1(\Sigma_g)$. Hence if we define $E^k(g)$ by the sequence
\[
\cdots \Lambda^{k-2}H_1(\Sigma_g)\stackrel{\omega\wedge\cdot}{\to} \Lambda^kH_1(\Sigma_g)\to E^k(g) \to \Lambda^{k-1}H_1(\Sigma_g)\stackrel{\omega\wedge\cdot}{\to} \Lambda^{k+1}H_1(\Sigma_g)\cdots
\]
then $E^k(g) \cong E_{2g+1-k}(g)$. Now, $\Lambda^kH_1(\Sigma_g)\cong H^k(T^{2g})$, where $T^{2g}$ is the Jacobian torus of $\Sigma_g$. Therefore, the sequence above can be identified with the Gysin sequence associated to the circle bundle $\E (g)$ over $T^{2g}$ having Euler class $\omega$, so that $E^*(g)=H^*(\E(g))$. 

The cohomology of $\E(g)$ was determined (in a different guise) by Lee and Packer \cite{leepacker}, using the Gysin sequence above together with combinatorial matrix theory. Taking the Poincar\'e dual of their result shows:
\[
E_k(g) \cong \left\{\begin{array}{ll} \zee^{\ts{{2g}\choose{k}} - {{2g}\choose{k-2}} }\oplus \ds \bigoplus_{j=2}^{\lfloor (k+1)/2\rfloor} \zee_j^{\ts{{2g}\choose{k-2j+1}} - {{2g}\choose{k-2j-1}}} & 0\leq k\leq g \\
\ds\bigoplus_{j = 0}^{\lfloor(2g+1-k)/2\rfloor} \zee_j^{\ts{{2g}\choose{k+2j-1}} - {{2g}\choose{k+2j+1}}} & g+1\leq k \leq 2g+1 \end{array}\right.
\]
where $\zee_0 = \zee$ and $\zee_1$ is the trivial group.

In particular it follows from this that for $k\leq g$
\[
\rk_{\mathbb Z} (E_k(g)) = \rk_{\mathbb Z} (E_{2g+1-k}(g)) =\ts {{2g}\choose{k}} - {{2g}\choose{k-2}},
\]
though it is possible to obtain the latter directly as well. It is an exercise based on this and \eqref{sigxs1homol} to see that
\[
\ts h(\Sigma_g\times S^1) = {{2g+1}\choose{g}}.
\]

We remark that the Floer homology $HF^\infty(\Sigma_g\times S^1,\s)$ for $c_1(\s) = 0$ was calculated in \cite{sigxs1} for coefficients in $\cee$ and in $\zee_2$; both these results are consistent with the hypothesis that $HF^\infty(\Sigma_g\times S^1,\s) \cong HC^\infty_*(\Sigma_g\times S^1)$, with coefficients in $\zee$.

More generally, consider a 3-manifold obtained as the mapping torus of a diffeomorphism $f: \Sigma_g\to \Sigma_g$. That is, $Y$ is constructed by gluing the boundaries of $\Sigma_g\times [0,1]$ via $f$. Then $H^1(Y) \cong \zee \oplus V$, where $\zee$ is generated by the Poincar\'e dual of $[\Sigma_g]$ and $V = \ker(1-f^*)$, $f^*$ denoting the action of $f$ on the first cohomology of $\Sigma_g$. It is not hard to see that the cup product form of $Y$ is given in this case by
\[
\mu_Y = s\wedge (\omega|_{V}),
\]
where $\omega$ is the intersection form on $\Sigma_g$ as before and $s$ is represented by a section of the obvious fibration $Y\to S^1$. Working over the rationals for simplicity, we can write $V = W\oplus V_0$, where $W$ is a maximal symplectic subspace of $V$ and $\omega|_{V_0} = 0$. Clearly this induces a decomposition
\[
C^\infty_*(Y) \cong C^\infty_*(\Sigma_w\times S^1)\otimes_{{\mathbb Z}[U,U^{-1}]} C^\infty_*(\#^{v_0}S^2\times S^1)
\]
of chain complexes, where $2w = \dim(W)$ and $v_0 = \dim(V_0)$. Applying the K\"unneth formula as in the proof of Theorem \ref{2ndthm} below and the results for $\#^n S^1\times S^2$ and $\Sigma_w\times S^1$ above, we infer
\[
\ts h(Y) = 2^{v_0}{{2w+1}\choose{w}}.
\]
One can see that this number is at least as large as the corresponding value of $h$ for a trivially fibered 3-manifold having the same first Betti number.

\section{Proofs}
\label{proofsec}

\subsection{Proof of Theorem \ref{mainthm}}\label{thm1proof}

Theorem \ref{mainthm} follows from a straightforward estimate of the size
of the homology of $(C^\infty_*,\partial)$. To state what we need explicitly, 
note that for each $i$, $0\leq i\leq b_1(Y)$, the differential restricts as a map
\[
\partial: \Lambda^iH^1(Y)\to \Lambda^{i-3}H^1(Y).
\]
Since the differential commutes with the action of $U$, we are 
reduced to considering three chain complexes $C_0$, $C_1$, and $C_2$, 
where
\begin{equation}
C_j = \bigoplus_{i\equiv j \mod 3} \Lambda^iH^1(Y).
\label{Cdef}
\end{equation}
We are interested in bounding the size of the homology of the $C_j$ from below; 
we do so by observing that the total rank of the homology of a 
chain complex must be at least the absolute value of the Euler 
characteristic of the complex.

\begin{prop} Fix an integer $b\geq 1$ and let $H$ denote a free Abelian group of rank $b$. Define graded groups $C_j$, $j = 0,1,2$ by the 
formula (\ref{Cdef}), where the grading on the factor $\Lambda^{3k + 
j}H$ is given by $k$. Then
\begin{itemize}
\item[a)] If $b$ is odd, then 
\[
|\chi(C_0)| + |\chi(C_1)| + |\chi(C_2)| = 2\cdot 3^{\frac{b-1}{2}}.
\]
\item[b)] If $b$ is even, then 
\[
|\chi(C_0)| + |\chi(C_1)| + |\chi(C_2)| = 4\cdot 3^{\frac{b}{2} - 1}.
\]
\end{itemize}
\label{mainprop}
\end{prop}

\begin{proof}[Proof of Theorem \ref{mainthm}] The upper bound on $h(Y)$ 
was proved in Lemma \ref{ublemma}. To 
obtain a lower bound, note that it suffices to 
consider two adjacent values of $k$, say $k = 0$ and $k = 1$.
It is easy to see that
\[
HC^\infty_0\oplus HC^\infty_1 \cong H_*(C_0)\oplus H_*(C_1)\oplus H_*(C_2).
\]
Since $\chi(H_*(C_j)) = \chi(C_j)$, the proposition, together with the 
obvious bound 
$\rk(H_*(C_j))\geq |\chi(H_*(C_j))|$, gives
\[
\rk(HC^\infty_0 \oplus HC^\infty_1) \geq |\chi(C_0)| + |\chi(C_1)| + 
|\chi(C_2)| = 2 L(b)
\]
where $b = b_1(Y)$, and Theorem \ref{mainthm} follows.\end{proof}

The proof of Proposition \ref{mainprop} is an exercise in the binomial 
theorem. To begin with, note that for $j = 0,1,2$,
\[
\chi(C_j) = \sum_k (-1)^k{{b}\choose{3k+j}}.
\]
To facilitate the discussion below, we denote the above sum by $S(b, 
j)$. Now, if $\xi$ satisfies $\xi^3 = 1$ then the binomial theorem gives
\begin{equation}
(1-\xi)^b = S(b,0) - S(b, 1)\xi + S(b,2)\xi^2.
\label{binthm}
\end{equation}
In particular, taking $\xi = 1$ we have
\begin{equation}
S(b,0) - S(b,1) + S(b,2) = 0.
\label{basicfact}
\end{equation}

Now we take $\xi = e^{2\pi i/3}$ and apply (\ref{binthm}) to the identity
\[
(1-\xi)^b = (1-\xi)(1-\xi)^{b-1}.
\]
Using $1-\xi = \frac{3}{2} - \frac{i\sqrt{3}}{2}$ and equating real 
and imaginary parts gives the relations
\begin{eqnarray*}
2S(b,0) + S(b,1) - S(b,2) &=& 3S(b-1,0) - 3S(b-1,2)\\
S(b,1) + S(b,2) &=& S(b-1,0) + 2 S(b-1,1) + S(b-1,2).
\end{eqnarray*}
Together with (\ref{basicfact}), this leads quickly to the recursion 
relations
\begin{eqnarray}
S(b,0)&=& S(b-1,0) - S(b-1,2)\label{recrel1} \\
S(b,1) &=& S(b-1,0) + S(b-1,1) \\
S(b,2) &=& S(b-1,1) + S(b-1,2).\label{recrel3}
\end{eqnarray}
Thus the values of $S(b,j)$ for $b$ even are determined by those for 
$b$ odd. We focus on the latter case.

First, an easy exercise using the symmetry ${{b}\choose{n}} = 
{{b}\choose{b-n}}$ shows:
\begin{eqnarray}
S(6n+1, 2) = 0 & \mbox{and} & S(6n+1, 0) = S(6n+1,1)\label{mod1rel}\\
S(6n+3, 0) = 0 & \mbox{and} & S(6n+3, 1) = S(6n+3, 2)\\
S(6n+5, 1) = 0 & \mbox{and} & S(6n+5, 0) = -S(6n+5, 2).
\end{eqnarray}

Next, we obtain a ``2-level'' recursion formula by applying 
(\ref{binthm}) to the identity
\[
(1-\xi)^b = (1-\xi)^2(1-\xi)^{b-2}
\]
using $(1-\xi)^2 = 3(\frac{1}{2} - \frac{i\sqrt{3}}{2})$ and equating 
real and imaginary parts as before. We get
\begin{eqnarray*}
\textstyle S(b,0) + \frac{1}{2}S(b,1) - \frac{1}{2}S(b,2) &=& 
\textstyle 3(\frac{1}{2}S(b-2,0) - \frac{1}{2}S(b-2,1) - S(b-2,2))\\
S(b,1) + S(b,2) &=& 3(S(b-2,0) + S(b-2,1)),
\end{eqnarray*}
which can be rearranged to give
\begin{eqnarray*}
S(b,0) + S(b,1) &=& 3(S(b-2,0) - S(b-2,2))\\
S(b,0) - S(b,2) &=& -3(S(b-2,1) + S(b-2,2))\\
S(b,1) + S(b,2) &=& 3(S(b-2,0) + S(b-2,1)).
\end{eqnarray*}
Substituting these equations into each other, we obtain
\begin{eqnarray*}
S(b,0) + S(b,1) &=& -3^3(S(b-6,0)+ S(b-6, 1))\\
S(b,0) - S(b,2) &=& -3^3(S(b-6,0) - S(b-6, 2))\\
S(b,1) + S(b,2) &=& -3^3(S(b-6, 1)+ S(b-6,2)).
\end{eqnarray*}
Now suppose that $b\equiv 1\mod 6$. Then according to (\ref{mod1rel}) 
we have $S(b,2) = 0$, so the third of the above 
equations shows that in this case 
$S(b,1)$ satisfies
\[
S(b,1) = -3^3 S(b-6, 1).
\]
Since $S(1,1)$ is obviously 1, we get that if $b \equiv 1 \mod 6$ 
then $S(b,1) = \pm 3^{(b-1)/2}$, where the sign depends on the value 
of $b$ modulo 12. Similar reasoning for the other odd values of $b$ 
modulo 6 shows:
\begin{eqnarray*}
\mbox{If $b\equiv 1\mod 6$}&:& S(b,2) = 0 \mbox{ and } S(b, 0) = 
S(b,1) = \pm 3^{\frac{b-1}{2}}\\
\mbox{If $b\equiv 3\mod 6$}&:& S(b,0) = 0 \mbox{ and } S(b,1)= 
S(b,2) = \pm 3^{\frac{b-1}{2}}\\
\mbox{If $b\equiv 5\mod 6$}&:& S(b,1) = 0 \mbox{ and } S(b,0) = 
-S(b,2) = \pm 3^{\frac{b-1}{2}}.
\end{eqnarray*}
In particular since in the notation of the proposition $S(b, j) = 
\chi(C_j)$, we have proved the proposition for the case of $b$ odd. 
The even case follows from this and the recursion relations 
(\ref{recrel1})--(\ref{recrel3}); for example, if 
$b\equiv 0\mod 6$ then
\[
S(b,0) = 2\cdot 3^{\frac{b}{2}-1} \quad\mbox{and}\quad S(b,1) = -S(b,2) = \pm 
3^{\frac{b}{2} - 1}.
\]
This completes the proof of Proposition \ref{mainprop}.

\subsection{Behavior under Connected Sum}\label{connsumproof}

\begin{proof}[Proof of Theorem \ref{2ndthm}] If $Y = Y_1\# Y_2$ is a connected sum, then we have a decomposition $H^1(Y) = H^1(Y_1)\oplus H^1(Y_2)$ and therefore
\[
\Lambda^*H^1(Y) \cong \Lambda^*H^1(Y_1)\otimes \Lambda^*H^1(Y_2).
\]
Under this decomposition, the cup-product form $\mu_Y$ satisfies
\[
\mu_Y = \mu_{Y_1}\otimes 1 + 1\otimes \mu_{Y_2},
\]
and since contraction is a derivation, $\partial_Y = \partial_{Y_1}\otimes 1 + (-1)^p 1\otimes \partial_{Y_2}$ on $C^\infty_p(Y_1) \otimes C^\infty_q(Y_2)$. Therefore we have a decomposition of chain complexes
\[
C^\infty_*(Y) \cong C^\infty_*(Y_1)\otimes_{{\mathbb Z}[U,U^{-1}]} C^\infty_*(Y_2).
\]

Working with coefficients in a field $\eff$, it follows that 
\begin{eqnarray*}
HC^\infty_0(Y_1\# Y_2) &\cong& (HC^\infty_0(Y_1)\otimes_{\mathbb F} HC^\infty_0(Y_2))\\ && \hspace*{1cm} \oplus(HC^\infty_1(Y_1)\otimes_{\mathbb F} HC^\infty_1(Y_2))
\end{eqnarray*}
and
\begin{eqnarray*}
HC^\infty_1(Y_1\# Y_2) &\cong& (HC^\infty_0(Y_1)\otimes_{\mathbb F} HC^\infty_1(Y_2))\\ && \hspace*{1cm}\oplus(HC^\infty_0(Y_1)\otimes_{\mathbb F} HC^\infty_1(Y_2)).
\end{eqnarray*}

In the notation of section \ref{prelimsec}, this gives
\begin{eqnarray*}
h_{ev}(Y_1\# Y_2) &=& h_{ev}(Y_1)h_{ev}(Y_2) + h_{odd}(Y_1)h_{odd}(Y_2)\\
h_{odd}(Y_1\# Y_2) &=& h_{ev}(Y_1)h_{odd}(Y_2) + h_{odd}(Y_1)h_{ev}(Y_2).
\end{eqnarray*}
Adding these equations yields
\[
h_{ev}(Y_1 \# Y_2) + h_{odd}(Y_1\# Y_2) = (h_{ev}(Y_1) + h_{odd}(Y_1))(h_{ev}(Y_2) + h_{odd}(Y_2)),
\]
and since $h = \frac{1}{2}(h_{ev} + h_{odd})$ this is
\[
2h(Y_1 \# Y_2) = 4h(Y_1)h(Y_2).
\]
\end{proof}

As a simple illustration of this result, we show that if $Y$ is a 3-manifold with $b_1(Y) = 5$ that is not rationally irreducible, then $h(Y)$ can only be 12 or 16. Indeed, suppose $Y$ decomposes as $Y = Y'\#Y''$ and assume first that $b_1(Y') = 2$ and $b_1(Y'') = 3$. Then $h(Y') = 2$, and $h(Y'')$ is either $3$ or $4$. These two cases give $h(Y) = 12$ or $16$, according to Theorem \ref{2ndthm}. The other possibility is that $b_1(Y') = 1$ and $b_1(Y'') = 4$: here $h(Y') = 1$ and $h(Y'')$ is $6$, $7$, or $8$. The case $h(Y'') = 7$ is ruled out by Lemma \ref{constrlemma}, and the other two cases again give $h(Y) = 12$ and $h(Y) = 16$.

\begin{proof}[Proof of Theorem \ref{irredthm}] We are given 3-manifolds $Y_1$ and
$Y_2$ with nonvanishing first Betti numbers $x = b_1(Y_1)$ and $y =
b_1(Y_2)$ that are not both odd. First suppose that $x$ and $y$ are
of opposite parity: say $x$ is odd and $y$ is even. Then
\[
h(Y_1\#Y_2) = 2h(Y_1)h(Y_2) \geq 2(3^{\frac{x-1}{2}}\cdot 2\cdot 
3^{\frac{y}{2} - 1}) = 4\cdot 3^{\frac{x+y-1}{2} - 1} = \frac{4}{3}L(b),
\]
where $b = x+y = b_1(Y_1\#Y_2)$ is odd.

Similarly, if both $x= b_1(Y_1)$ and $y = b_1(Y_2)$ are even, we
have
\[
h(Y_1\#Y_2) = 2h(Y_1)h(Y_2) \geq 8\cdot 3^{\frac{x}{2}-1}\cdot 
3^{\frac{y}{2}-1} = \frac{8}{3}\cdot 3^{\frac{x+y}{2} - 1} = 
\frac{4}{3}L(b).
\]
\end{proof}

\section{Relation to Floer homology}\label{floersec}

We outline some results of Ozsv\'ath and Szab\'o relating to the structure of $HF^\infty(Y,\s)$ for $\s$ a torsion \spinc structure, and describe the relationship to $HC^\infty_*(Y)$.

Recall that there is a version of Heegaard Floer homology with ``universal'' coefficients in the ring $R_Y = \zee[H^1(Y;\zee)]$. We have the following general result of \OnS:

\begin{thm}[Theorem 10.12 of \cite{OS2}]\label{twistcoeffthm} If $(Y,\s)$ is a closed \spinc 3-manifold and $c_1(\s)$ is torsion, then there is an isomorphism
\[
HF^\infty(Y,\s; R_Y) \cong \zee[U,U^{-1}]
\]
of $R_Y$-modules, where elements of $H^1(Y;\zee)$ act as the identity on $\zee[U,U^{-1}]$. 
\end{thm}

There is an additional structure on Heegaard Floer homology in a {\spinc} structure with torsion Chern class: a grading that takes values in the rational numbers, constructed by Ozsv\'ath and Szab\'o in \cite{OS3}. In the case of $HF^\infty$, the theorem above provides an integral grading, natural up to an integer shift, with respect to which nonvanishing homogeneous parts lie in even degrees. The relationship between these two gradings is determined by the \spinc structure $\s$ via the formula for the shift in grading induced by cobordisms \cite{OS3}: explicitly, write the \spinc 3-manifold $(Y,\s)$ as the \spinc boundary of a 4-manifold $(Z,\R)$. Then the rational-valued grading on the Heegaard Floer homology of $(Y,\s)$ takes values in $\zee + r$, where $r\in \cue$ is given by
\[
r = \frac{1}{4}(c_1^2(\R) - 3\sigma(Z) - 2e(Z)) + \frac{1}{2}.
\]
Here $\sigma(Z)$ is the signature of the intersection form of $Z$ and $e(Z)$ is the Euler characteristic. 

For our purposes the rational-valued grading is unimportant, so we simply impose an integer grading on $HF^\infty$ in such a way that the universally-twisted homology is supported in even degrees as above. This choice induces an integer grading on $HF^\infty(Y,\s,M)$ for any coefficient module $M$---in particular for $M = \zee$---that is well-defined up to a shift by an even integer.

The $\zee$-coefficient Floer homology can be calculated from the universal version by making use of a change-of-coefficients spectral sequence. For the benefit of readers not familiar with this construction, we outline it here in general before considering the case at hand. 

Suppose $(C_*, d)$ is a free chain complex of modules over a commutative ring $R$, and $M$ is an $R$-module: we wish to compute the homology of the complex $C_*\otimes_R M$. Under reasonable circumstances we may assume $M$ has free resolution
\[
\cdots \stackrel{\delta_2}{\longto} P_2\stackrel{\delta_1}{\longto} 
P_1\stackrel{\delta_0}{\longto} P_0\to M\to 0,
\]
i.e., each $P_j$ is a free $R$-module, and the sequence above is exact (a projective resolution would also suffice). We can form the double complex $E^0_{i,j} = C_i \otimes_R P_j$, with the obvious pair of differentials $d: E^0_{i,j}\to E^0_{i-1,j}$ (the ``horizontal'' differential) 
and $\delta: E^0_{i,j}\to E^0_{i,j-1}$ (the ``vertical'' differential) and ``total complex'' 
$(\tot(E)_*,D)$ where $\tot(E)_k = \bigoplus_{i+j=k}E^0_{i,j}$ and $D 
= d+ (-1)^j\delta$ on $E_{i,j}^0$. 

The double complex gives rise to a spectral sequence $(E^n, d_n)$ converging to the associated graded module $Gr(H_*(\tot(E)))$ determined by one of two natural filtrations on $\tot(E)$: the 
``horizontal'' filtration $\cdots\subset F^h_0\subset F^h_1\subset \cdots$ where $F^h_j = 
\bigoplus_{j'\leq j} E^0_{i,j'}$, and the ``vertical'' filtration 
$\cdots \subset F_0^v\subset F_1^v\subset\cdots$, $F_i^v = \bigoplus_{i'\leq i} E^0_{i',j}$. In fact, each of these filtrations gives rise to its own spectral sequence, but as we shall see one of these is essentially trivial: hence we will suppress the filtration from the notation for the sequence $(E^n, d_n)$.

The spectral sequence is constructed as follows. Choose one of the filtrations $F_j^h$ or $F_i^v$ and take the homology of $E^0$ with respect to the corresponding differential. The result is called $E^1_{i,j}$, and it inherits a differential from the original complex, corresponding to the differential in the direction not used at first. For 
example, with the double complex $E^0_{i,j}= C_i\otimes P_j$, consider 
the homology with respect to the differential $\delta$---that 
is, the ``vertical'' homology.  Since the $C_i$ are free modules, the 
result is $E^1_{i,j} = C_i\otimes H_j(P_*)$.  Of course, by 
construction $H_j(P_*) = 0$ except when $j=0$, and $H_0(P_*) = M$.  
Thus starting with the vertical differential gives a spectral sequence with 
\[
E^1_{i,j} = \left\{ \begin{array}{ll} 0 & \mbox{if $j>0$} \\
C_i\otimes M & \mbox{if $j=0$}.\end{array}\right.
\]
Taking the homology with respect to the remaining differential, the 
``horizontal'' one, gives the $E^2$ term of this sequence, which is obviously
\[
E^2_{i,j} = \left\{ \begin{array}{ll} 0 & \mbox{if $j>0$}\\
H_i(C_*\otimes M) & \mbox{if $j=0$}.\end{array}\right.
\]
The general machinery of spectral sequences would now give a 
differential $d_2: E^2_{i,j} \to E^2_{i-2, j+1}$ whose homology is the 
$E^3$-term, but from the structure of $E^2$ above, this differential 
must be trivial (and the same for all subsequent differentials).  The spectral sequence corresponding to this 
filtration therefore collapses at the $E^2$ term, and we infer 
that $H_*(\tot(E)) = H_*(C_*\otimes M)$ is the homology we wish to calculate.

To calculate this homology, we turn to the other filtration---that is, 
we use the other differential first.  Returning to $E^0$, we take the 
homology with respect to the horizontal differential, the differential 
of $C_*$.  Again, since the $P_j$ are free, we get that 
\[
E^1_{i,j} = H_i(C_*)\otimes_R P_j.
\]
When we take the next homology to get $E^2$, that is the homology in 
the vertical direction, coming from the $\delta_j$, we can no longer 
commute homology and tensor product (since the $H_i(C_*)$ need not be free $R$-modules): by definition, the vertical homology is
\[
E^2_{i,j} = \Tor^R_j(H_i(C_*),M).
\]
We have recovered the following standard fact.
\begin{prop} Given a free chain complex of $R$-modules $C_*$ and 
another module $M$, there is a ``universal coefficients spectral sequence'' 
converging to the homology $H_*(C_*\otimes_R M)$ whose $E^2$ term is 
given by $E^2_{i,j} = \Tor_j^R(H_i(C_*),M)$, with differential $d_2: 
E^2_{i,j} \to E^2_{i+1, j-2}$.
\end{prop}

Returning now to Heegaard Floer homology, Theorem \ref{twistcoeffthm} states that the Floer homology $HF^\infty(Y,\s, R_Y)$ for $\s$ a \spinc structure with $c_1(\s)$ torsion is equal to $0$ or $\zee$ in alternating degrees. Hence the $E^2$ term of the universal coefficient spectral sequence (taking $M = \zee$) has zeros in each odd column, while the remaining columns are given by
\[
E^2_{i,*} = \Tor^{R_Y}_*(\zee,\zee)
\]
for each even $i$. It is a standard fact that if $G$ is an abelian group then 
\[
\Tor^{{\mathbb Z}[G]}_*(\zee,\zee) = H_*(G;\zee) = H_*(K(G,1);\zee).
\]
In our case, $G = H^1(Y;\zee)$ is free abelian of rank $b_1(Y)$, and $K(G,1)$ is the torus $T^b$ of dimension $b = b_1(Y)$. Thus we have a natural identification $\Tor^{R_Y}_*(\zee,\zee) = H_*(T^b;\zee) \cong \Lambda^*H^1(Y;\zee)$, and the $E^2$ term of the universal coefficients spectral sequence can be written 
\[
\mbox{Tor}_*^{R_Y}(HF^\infty_*(Y,\s;R_Y), \zee) \cong \Lambda^*H^1(Y;\zee)\otimes \zee[U,U^{-1}]\cong C^\infty_*(Y).
\] 
Observe that according to the proposition above, the $d_2$ differential maps one column to the right and two rows down. Since every other column in the $E^2$ complex vanishes we infer that $d_2 = 0$, so that the above is also the $E^3$ term of the spectral sequence. The $d_3$ differential restricts as a map $d_3:\Lambda^kH^1(Y)\otimes U^n\to \Lambda^{k-3}H^1(Y)\otimes U^{n-1}$. {\OnS} conjecture in \cite{OSplumb} that $(E^3, d_3)$ is the complex $(C^\infty_*(Y),\partial)$ considered in this article.

Furthermore, it is conjectured in \cite{OSplumb} that all subsequent differentials in the spectral sequence vanish, so that the $E^\infty$ term in the universal coefficients spectral sequence is our $HC^\infty_*(Y)$. 

While we do not address these conjectures here, we observe that since the arguments in the proof of Theorem \ref{mainthm} do not make use of the differential on $C^\infty_*(Y)$ but only the ranks of the chain groups, the bounds obtained there apply to the rank in each degree of $HF^\infty(Y,\s;\zee)$ provided that the universal coefficients sequence collapses after the $E_3$ stage (see also the remarks on monopole Floer homology at the end of the introduction). Let us say that $(Y,\s)$ is {\it regular} if $c_1(\s)$ is a torsion class and all differentials $d_r$, $r\geq 4$, in that spectral sequence vanish.

\begin{prop} If $(Y,\s)$ is a regular \spinc 3-manifold then the rank of $HF^\infty_k(Y,\s)$ satisfies
\[
L(b)\leq \rk \,HF^\infty_k(Y;\s)\leq 2^{b-1},
\]
where $b = b_1(Y)$. 
\end{prop}\hfill$\Box$

The even differentials $d_{2r}$ in the sequence we are considering vanish for dimensional reasons; thus the first differential past $d_3$ that may be nontrivial is $d_5: \Lambda^kH^1(Y)\otimes U^n\to \Lambda^{k-5}H^1(Y)\otimes U^{n-2}$. Hence:

\begin{prop} If $b_1(Y)\leq 4$ then $(Y,\s)$ is regular for any torsion {\spinc} structure $\s$.
\end{prop}\hfill$\Box$

This proves Corollary \ref{minhcor}.


\begin{thebibliography}{99}
\bibitem{hedden} Matthew Hedden, ``On knot Floer homology and cabling,'' {\it Algebr. Geom. Topol.} {\bf 5} (2005), 1197--1222.
\bibitem{sigxs1} Stanislav Jabuka and Thomas Mark, ``On the Heegaard Floer homology of a surface times a circle,'' arXiv:math.GT/0502328.
\bibitem{KM} Peter Kronheimer and Tomasz Mrowka, {\it Monopoles and three-manifolds,} book in preparation.
\bibitem{leepacker} Soo Teck Lee and Judith A. Packer, ``The cohomology of the integer Heisenberg groups,'' {\it J. Algebra} {\bf 184} (1996), no. 1, 230--250.
\bibitem{stipsicz} Paolo Lisca and Andr\'as Stipsicz, ``Ozsv\'ath-Szab\'o invariants and tight contact three-manifolds. I.'' {\it Geom. Topol.} {\bf 8} (2004), 925--945.
\bibitem{nemethi} Andr\'as Nemethi, ``On the Ozsv\'ath-Szab\'o invariant of negative definite plumbed 3-manifolds,'' {\it Geom. Topol.} {\bf 9} (2005), 991--1042.
\bibitem{OS4} Peter Ozsv\'ath and Zolt\'an Szab\'o, ``Absolutely graded Floer homologies and intersection forms for four-manifolds with boundary,'' {\it Adv. Math.} {\bf 173} (2003), no. 2, 179--261.
\bibitem{OSplumb} Peter Ozsv\'ath and Zolt\'an Szab\'o, ``On the Floer homology of plumbed three-manifolds,'' {\it Geom. Topol.} {\bf 7} (2003), 185--224.
\bibitem{OS1} Peter Ozsv\'ath and Zolt\'an Szab\'o, ``Holomorphic disks and topological invariants for closed three-manifolds,'' {\it Ann. of Math.} {\bf 159} (2004), no. 3, 1027--1158.
\bibitem{OS2} Peter Ozsv\'ath and Zolt\'an Szab\'o, ``Holomorphic disks and three-manifold invariants: properties and applications,'' {\it Ann. of Math.} {\bf 159} (2004), no. 3, 1159--1245.
\bibitem{OS3} Peter Ozsv\'ath and Zoltan Szab\'o, ``Holomorphic triangles and invariants for smooth four-manifolds,'' {\it Adv. Math.} {\bf 202} (2006), no. 2, 326--400.
\bibitem{sullivan} Dennis Sullivan, ``On the intersection ring of compact three manifolds,'' {\it Topology} {\bf 14} (1975) 275--277.
\end{thebibliography}
\end{document}